\documentclass[11pt, a4paper, twoside]{article}
\usepackage{amsfonts}
\usepackage{color}
\usepackage{mathrsfs}
\usepackage{amsmath}
\usepackage{amssymb}
\usepackage{fancyhdr}
\usepackage{hyperref}

\numberwithin{equation}{section}

\allowdisplaybreaks

\begin{document}

\fancyhf{}

\fancyhead[EC]{M. Song}

\fancyhead[EL]{\thepage}

\fancyhead[OC]{Strichartz inequalities for the Schr\"{o}dinger equation}

\fancyhead[OR]{\thepage}

\renewcommand{\headrulewidth}{0pt}
\renewcommand{\thefootnote}{\fnsymbol {footnote}}

\title{\textbf{Strichartz inequalities for the Schr\"{o}dinger equation with the full Laplacian on H-type groups}}
\footnotetext {{}\emph{Key words and phrases}: Strichartz inequalities, Schr\"{o}dinger equation, full Laplacian, H-type groups.}
\setcounter{footnote}{0}
\author{Manli Song}
\date{}
\maketitle

\begin{abstract}
In this paper, we prove the dispersive estimates and Strichartz inequalities for the solution of the Schr\"{o}dinger equation related to the full Laplacian on H-type groups. This
extends the results obtained by G. Furioli and A. Veneruso [Studia Math., 2004, 160: 157-178] on the Heisenberg group.
\end{abstract}

\newtheorem{theorem}{Theorem}[section]
\newtheorem{preliminaries}{Preliminaries}[section]
\newtheorem{definition}{Difinition}[section]
\newtheorem{main result}{Main Result}[section]
\newtheorem{lemma}{Lemma}[section]
\newtheorem{proposition}{Proposition}[section]
\newtheorem{corollary}{Corollary}[section]
\newtheorem{remark}{Remark}[section]

\section[Introduction]{Introduction}
The aim of the paper is to study the Stricharz estimates for the following Cauchy problem of the Schr\"{o}dinger equation related to the full Laplacian on H-type groups $G$:
\begin{equation}\label{schrodinger}
\left\{
\begin{array}{ll}
&\partial_tu-i\mathcal{L}u=f \in L^1((0, T),L^2)\\
&u|_{t=0}=u_0 \in \dot{B}^1_{2,2},
\end{array}\right.
\end{equation}
where $G$ is the H-type group with topological dimension $n$ and homogeneous dimension $N$, $\mathcal{L}$ is its full Laplacian and the Besov spaces $\dot{B}^\rho_{q,r}$ are defined by a Littlewood-Paley decomposition related to $\mathcal{L}$.

Strichartz estimates are very useful in the study of nonlinear partial differential equations. These estimates in the Euclidean setting have been proved for many dispersive equations, such as wave equation and Schr\"{o}dinger equation (see \cite{GV, KT, Str}). To obtain Strichartz estimates, it involves basically two types of ingredients. The first one consists in estimating the decay in time on the evolution group associated with the free equation (i.e. $f=0$). The second one consists of abstract arguments, which are mainly duality arguments.

Many authors have also been interested in adapting the well known Strichartz estimates from the Euclidean setting to a more abstract setting, such as the Heisenberg group and the H-type groups. In 2000, H. Bahouri, P. G\'{e}rard and C.-J. Xu \cite{BGX} discussed the Strichartz estimates with the sublaplacian on the Heisenberg group, by means of Besov spaces defined by a Littlewood-Paley decomposition related to the spectral of the sublaplacian. In their work, they showed such estimates existed for the wave equation while failed for the Schr\"{o}dinger equation. To avoid the particular behavior of the Schr\"{o}dinger operator on the Heisenberg group, the sublaplacian was been replaced by the full Laplacian (see \cite{FMV1, FV}). Later, Strichartz estimates were addressed on general H-type groups (see \cite{DH, Song, SZ}), but they only considered Besov spaces related to the sublaplacian. In a recent paper \cite{LS}, H. Liu and M. Song have proved the Strichartz inequalities for the wave equation with the full Laplacian on H-type groups. Thus, it is natural to wonder whether such estimates also remain true for the corresponding Schr\"{o}dinger equation.

Our purpose is to show that the Schr\"{o}dinger equation related to the full Laplacian on H-type groups is also dispersive. In comparison with \cite{DH}, the full Laplacian does not have the homogeneous properties, which involves some technical difficulties. Furthermore, in those groups, only the Heisenberg group has a one dimensional center. Let $p$ be the dimension of the center on H-type groups. In this paper, we only consider those groups with $p>1$, which makes the issue become very complicated.

It is well-known that the solution of the non-homogeneous equation \eqref{schrodinger} is given by the sum $u=v+w$ where
\begin{equation*}
v(t)=e^{it\mathcal{L}}u_0
\end{equation*}
is the solution of \eqref{schrodinger} with $f=0$ and
\begin{equation*}
w(t)=\int_0^te^{i(t-\tau)\mathcal{L}}f(\tau)\,d\tau
\end{equation*}
is the solution of \eqref{schrodinger} with $u_0=0$.

We can now state the main results of the paper. We first give the sharp dispersive estimate on the free solution.
\begin{theorem} \label{dispersive inequ} \indent If $v$ is the free solution of the Schr\"{o}dinger equation \eqref{schrodinger}, then
\begin{equation*}
||v(t)||_{L^\infty(G)}\leq C|t|^{-p/2}||u_0||_{\dot{B}^{n-1}_{1,1}},
\end{equation*}
and the result is sharp in time.
\end{theorem}

In comparison with the results for the Schr\"{o}dinger equation in \cite{FV} and \cite{DH}, we have obtained an improvement on the time decay, respectively due to the replacement of the full Laplacian and the bigger size of the H-type group center.

We also get a very useful estimate on the solution of the Schr\"{o}dinger equation.
\begin{theorem}\label{solution}\indent For $i=1,2$, let $q_i, r_i\in[2,\infty]$ and $\rho_i\in\mathbb{R}$ such that
\begin{align*}
&a)\,\,\frac{2}{q_i}=p(\frac{1}{2}-\frac{1}{r_i}); \\
&b)\,\,\rho_i=-(n-1)(\frac{1}{2}-\frac{1}{r_i}),
\end{align*}
except for $(q_i, r_i, p)=(2,\infty,2)$. Let $q_i'$, $r_i'$ denote the conjugate exponent of $q_i$ and $r_i$. The solution of the Cauchy problem \eqref{schrodinger} $u$ satisfies the estimate
\begin{equation*}
||u||_{L^{q_1}((0,T),\dot{B}^{\rho_1}_{r_1,2})}\leq C\big(||u_0||_{L^2(G)}+||f||_{L^{q_2'}((0,T),\dot{B}^{-\rho_2}_{r_2',2})}\big),
\end{equation*}
where the constant $C>0$ does not depend on $u_0$, $f$ or $T$.
\end{theorem}

\section[H-type groups and spherical Fourier transform]{H-type groups and spherical Fourier transform}
\begin{bf} {2.1. H-type groups.}\end{bf} \indent Let $\mathfrak{g}$ be a two step nilpotent Lie algebra endowed with an inner product $\langle \cdot,\cdot \rangle$. Its center is denoted by $\mathfrak{z}$. $\mathfrak{g}$ is said to be of H-type if $[\mathfrak{z}^{\bot},\mathfrak{z}^{\bot}]=\mathfrak{z}$ and for every $s \in \mathfrak{z}$, the map $J_s: \mathfrak{z}^{\bot} \rightarrow \mathfrak{z}^{\bot}$ defined by
\begin{equation*}
\langle J_s u, w \rangle:=\langle s, [u,w] \rangle, \quad \forall u, w \in \mathfrak{z}^{\bot},
\end{equation*}
is an orthogonal map whenever $|s|=1$. An H-type group is a connected and simply connected Lie group $G$ whose Lie algebra is of H-type.

Given $0 \neq a \in \mathfrak{z}^*$, the dual of $\mathfrak{z}$, we can define a skew-symmetric mapping $B(a)$ on $\mathfrak{z}^{\bot}$ by
\begin{equation*}
\langle B(a)u,w \rangle =a([u,w]), \quad \forall u,w \in \mathfrak{z}^{\bot}.\\
\end{equation*}
We denote by $z_a$ the element of $\mathfrak{z}$ determined by
\begin{equation*}
\langle B(a)u,w \rangle =a([u,w])=\langle J_{z_a} u,w \rangle. \\
\end{equation*}
Since $B(a)$ is skew symmetric and non-degenerate, the dimension of $\mathfrak{z}^{\bot}$ is even, i.e., $dim\mathfrak{z}^{\bot}=2d$.

We can choose an orthonormal basis
\begin{equation*}
\{E_1(a),E_2(a),\cdots,E_d(a),\overline{E}_1(a),\overline{E}_2(a),\cdots,\overline{E}_d(a)\},
\end{equation*}
of $\mathfrak{z}^{\bot}$ such that
\begin{equation*}
B(a)E_i(a)=|z_a|J_{\frac{z_a}{|z_a|}}E_i(a)=|a|\overline{E}_i(a)
\end{equation*}
and
\begin{equation*}
B(a)\overline{E}_i(a)=-|a|E_i(a).
\end{equation*}

We set $p=dim \mathfrak{z}$. Throughout this paper we assume that $p>1$. We can choose an orthonormal basis $\{\epsilon_1,\epsilon_2,\cdots,\epsilon_p \}$ of $\mathfrak{z}$ such that $a(\epsilon_1)=|a|,a(\epsilon_j)=0,j=2,3,\cdots,p$. Then we can denote the element of $\mathfrak{g}$ by
\begin{equation*}
(z,t)=(x,y,t)=\underset{i=1}{\overset{d}{\sum}}(x_i E_i+y_i \overline{E}_i )+\underset{j=1}{\overset{p}{\sum}}s_j \epsilon_j.
\end{equation*}

We identify G with its Lie algebra $\mathfrak{g}$ by exponential map. The group law on H-type group $G$ has the form
\begin{equation}
(z,s)(z',s')=(z+z',s+s'+\frac{1}{2}[z,z']),  \label{equ:Law}
\end{equation}
where $[z,z']_j=\langle z,U^jz' \rangle$ for a suitable skew symmetric matrix $U^j,j=1,2,\cdots,p$.

\begin{theorem} \indent G is an H-type group with underlying manifold $\mathbb{R}^{2d+p}$, with the group law  $\eqref{equ:Law}$ and the matrix $U^j,j=1,2,\cdots, p$ satisfies the following conditions:\\
$(i)$ $U^j$ is a $2d \times 2d$ skew symmetric and orthogonal matrix, $j=1,2,\cdots, p$;\\
$(ii)$ $U^i U^j+U^j U^i=0,i,j=1,2,\cdots,p$ with $i \neq j$.
\end{theorem}

{\bf Proof.} See \cite{BU}.

\begin{remark} It is well know that H-type algebras are closely related to Clifford modules (see \cite{R}). H-type algebras can be classified by the standard theory of Clifford algebras. Specially, on H-type group $G$, there is a relation between the dimension of the center and its orthogonal complement space. That is $p+1\leq 2d$ (see \cite{A,KR}).
\end{remark}

\begin{remark}
We identify $G$ with $\mathbb{R}^{2d}\times\mathbb{R}^p$. We shall denote the topological dimension of $G$ by $n=2d+p$. Following Folland and Stein (see \cite{FS}), we will exploit the canonical homogeneous structure, given by the family of dilations$\{\delta_r\}_{r>0}$,
\begin{equation*}
\delta_r(z,s)=(rz,r^2s).
\end{equation*}
We then define the homogeneous dimension of $G$ by $N=2d+2p$.
\end{remark}

The left invariant vector fields which agree respectively with $\frac{\partial}{\partial x_j},\frac{\partial}{\partial y_j}$ at the origin are given by
\begin{equation*}
\begin{aligned}
X_j&=\frac{\partial}{\partial x_j}+\frac{1}{2}\underset{k=1}{\overset{p}{\sum}} \left( \underset{l=1}{\overset{2d}{\sum}}z_l U_{l,j}^k \right) \frac{\partial}{\partial s_k},\\
Y_j&=\frac{\partial}{\partial y_j}+\frac{1}{2}\underset{k=1}{\overset{p}{\sum}} \left( \underset{l=1}{\overset{2d}{\sum}}z_l U_{l,j+d}^k \right) \frac{\partial}{\partial s_k},\\
\end{aligned}
\end{equation*}
where $z_l=x_l,z_{l+d}=y_l,l=1,2,\cdots, d.$

The vector fields $S_k=\frac{\partial}{\partial s_k},k=1,2,\cdots,p$ correspond to the center of $G$. In terms of these vector fields we introduce the sublaplacian $\Delta$ and full Laplacian $\mathcal{L}$ respectively
\begin{equation}
\label{eq:Lap}
\begin{aligned}
\Delta&=-\underset{j=1}{\overset{d}{\sum}}(X_j^2 +Y_j^2)=\Delta_z +\frac{1}{4} |z|^2\mathcal{S}-\underset{k=1}{\overset{p}{\sum}}\langle z,U^k \nabla_z \rangle S_k,  \\
\mathcal{L}&=\Delta+\mathcal{S},  \\
\end{aligned}
\end{equation}
where
\begin{equation*}
\Delta_z=-\underset{j=1}{\overset{2d}{\sum}}\frac{\partial^2}{\partial z_j^2}, \mathcal{S}=-\underset{k=1}{\overset{p}{\sum}}\frac{\partial^2}{\partial s_k^2},\nabla_z=(\frac{\partial}{\partial z_1},\frac{\partial}{\partial z_2},\cdots,\frac{\partial}{\partial z_{2d}})^t.\\
\end{equation*}

\begin{bf} {2.2. Spherical Fourier transform.}\end{bf} The spherical functions associated to the Gelfand pair $(G, O(d))$ (we identify $O(d)$ with $O(d)\otimes Id_p$) have been computed in \cite{DR} and \cite{Kor}. They involve, as on the Heisenberg group, the Laguerre functions
\begin{equation*}
\mathfrak{L}_m^{(\alpha)}(\tau)=L_m^{(\alpha)}(\tau)e^{-\tau/2}, \quad\tau \in \mathbb{R}, m,\alpha \in \mathbb{N},
\end{equation*}
where $L_m^{(\alpha)}$ is the Laguerre polynomial of type $\alpha$ and degree $m$.

We say a function $f$ on $G$ is radial if the value of $f(z,s)$ depends only on $|z|$ and $s$. In the sequel, we denote by $\mathcal{S}_{rad}(G)$ and $L^q_{rad}(G)$,$1\leq q\leq \infty$, the spaces of radial functions in $\mathcal{S}(G)$ and $L^p(G)$, respectively. In particular, $L^1_{rad}(G)$ endowed with the convolution product
\begin{equation*}
f_1*f_2(g)=\int_Gf_1(gg'^{-1})f_2(g')\,dg',\quad g\in G
\end{equation*}
is a commutative algebra.

Let $f\in L^1_{rad}(G)$, and we define the spherical Fourier transform by
\begin{equation*}
\mathcal{F}(f)(\lambda,m)=\hat{f}(\lambda,m)=\left( \begin{array}{c} m+d-1\\m \end{array} \right)^{-1}
\int_{\mathbb{R}^{2d+p}}e^{i\lambda s} f(z,s)\mathfrak{L}_m^{(d-1)}(\frac{|\lambda|}{2}|z|^2)\,dzds.
\end{equation*}

A straightforward computation implies that $\mathcal{F}(f_1*f_2)=\mathcal{F}(f_1)\mathcal{F}(f_2)$. As in \cite{DH}, we also deduce the corresponding Fourier inversion formula.
\begin{proposition} \indent For all $f \in \mathcal{S}_{rad}(G)$ such that
\begin{equation*}
\underset{m\in\mathbb{N}}{\sum}\left( \begin{array}{c} m+d-1\\m \end{array} \right) \int_{\mathbb{R}^p} |\hat{f}(\lambda,m)||\lambda|^d d\lambda <\infty,
\end{equation*}
we have
\begin{equation}
f(z,s)=(\frac{1}{2\pi})^{d+p}\underset{m\in\mathbb{N}}{\sum}\int_{\mathbb{R}^p} e^{-i\lambda s} \hat{f}(\lambda,m) \mathfrak{L}_m^{(d-1)}(\frac{|\lambda|}{2}|z|^2)|\lambda|^d \,d\lambda, \label{Fourier-inversion}
\end{equation}
the sum being convergent in $L^{\infty}$ norm.
\end{proposition}

Moreover, if $f \in \mathcal{S}_{rad}(G)$, then $\mathcal{L}f\in\mathcal{S}_{rad}(G)$ and its spherical Fourier transform is given by
\begin{equation*}
\widehat{\mathcal{L}f}(\lambda,m)=((2m+d)|\lambda|+|\lambda|^2)\hat{f}(\lambda,m).
\end{equation*}
The full Laplacian $\mathcal{L}$ is a positive self-adjoint operator densely defined on $L^2(G)$. So by the spectral theorem, for any bounded Borel function $h$ on $\mathbb{R}$, we have
\begin{equation*}
\widehat{h(\mathcal{L})f}(\lambda,m)=h((2m+d)|\lambda|+|\lambda|^2)\hat{f}(\lambda,m).
\end{equation*}

\section[Homogeneous Besov spaces]{Homogeneous Besov spaces}
We shall recall the homogeneous Besov spaces given in \cite{LS}. Let $R$ be a non-negative, even function in $C_c^{\infty}(\mathbb{R})$ such that supp$R \subseteq \{\tau \in \mathbb{R}:\frac{1}{2}\leq |\tau|\leq4 \}$ and
\begin{equation*}
\underset{j\in \mathbb{Z}}{\sum}R(2^{-2j}\tau)=1, \quad\forall \tau \neq0.
\end{equation*}
For $j \in \mathbb{Z}$, we denote by $\psi_j$ the kernel of the operator $R(2^{-2j}\mathcal{L})$ and set $\Delta_j f=f*\psi_j$. As $R \in C_c^{\infty}(\mathbb{R})$ , A. Hulanicki \cite{Hul} proved that $\psi_j\in \mathcal{S}_{rad}(G)$ and
\begin{equation*}
\hat{\psi}_j(\lambda,m)=R(2^{-2j}((2m+d)|\lambda|+|\lambda|^2)).
\end{equation*}
By \cite{FMV2} (see Proposition 6), there exists $C>0$ such that
\begin{equation}\label{psi}
||\psi_j||_{L^1(G)}\leq C,\quad \forall j \in \mathbb{Z}.
\end{equation}
By standard arguments (see \cite{FMV2}, Proposition 9), we can deduce from \eqref{psi} that
\begin{equation}\label{sigma}
||\mathcal{L}^{\gamma/2}\Delta_j f||_{L^q(G)}\leq C2^{j\gamma}||\Delta_j f||_{L^q(G)}, \quad\gamma\in\mathbb{R}, j\in\mathbb{Z}, 1\leq q\leq\infty, f\in\mathcal{S}'(G),
\end{equation}
where both sides of \eqref{sigma} are allowed to be infinite.

By the spectral theorem, for any $f\in L^2(G)$, the following homogeneous Littlewood-Paley decomposition holds:
\begin{equation*}
f=\sum_{j\in\mathbb{Z}}\Delta_j f \quad \text{in $L^2(G)$}.
\end{equation*}
So
\begin{equation}\label{infty}
||f||_{L^\infty(G)}\leq\sum_{j\in\mathbb{Z}}||\Delta_j f||_{L^\infty(G)},\quad f\in L^2(G),
\end{equation}
where both sides of \eqref{infty} are allowed to be infinite.

Let $1\leq q,r \leq \infty, \rho <N/q$, we define the homogeneous Besov space $\dot{B}^\rho_{q,r}$ as the set of distributions $f \in \mathcal{S}^{'}(G)$ such that \begin{equation*}
||f||_{\dot{B}^\rho_{q,r}}=\left(\underset{j\in \mathbb{Z}}{\sum}2^{j\rho r}||\Delta_jf||_q^r\right)^{\frac{1}{r}}<\infty,
\end{equation*}
and $f=\underset{j\in \mathbb{Z}}{\sum}\Delta_jf$ in $\mathcal{S}'(G)$.

We collect all the properties we need about the spaces $\dot{B}^\rho_{q,r}$ in the following proposition.
\begin{proposition}Let $q,r\in [1,\infty]$ and $\rho<\frac{N}{q}$.\\
(i) The space $\dot{B}^\rho_{q,r}$ is a Banach space with the norm $||\cdot||_{\dot{B}^\rho_{q,r}}$;\\
(ii) the definition of $\dot{B}^\rho_{q,r}$ does not depend on the choice of the function $R$ in the Littlewood-Paley \\
\indent  decomposition;\\
(iii) for $-\frac{N}{q'}<\rho<\frac{N}{q}$, the dual space of $\dot{B}^\rho_{q,r}$ is $\dot{B}^{-\rho}_{q',r'}$;\\
(iv) for $\alpha\in[n,N]$, we have the continuous inclusion
\begin{equation*}
\dot{B}^{\rho_1}_{q_1,r}\subset\dot{B}^{\rho_2}_{q_2,r}, \quad \frac{1}{q_1}-\frac{\rho_1}{\alpha}=\frac{1}{q_2}-\frac{\rho_2}{\alpha}, \quad\rho_1\geq\rho_2;
\end{equation*}
(v) for all $q\in[2, \infty]$, we have the continuous inclusion $\dot{B}^0_{q,2}\subset L^q$;\\
(vi) $\dot{B}^0_{2,2}=L^2$;\\
(vii) for $\theta\in[0,1]$, we have
\begin{equation*}
[\dot{B}^{\rho_1}_{q_1,r_1}, \dot{B}^{\rho_2}_{q_2,r_2}]_\theta=\dot{B}^\rho_{q,r},
\end{equation*}
\indent with $\rho=(1-\theta)\rho_1+\theta\rho_2$, $\frac{1}{q}=\frac{1-\theta}{q_1}+\frac{\theta}{q_2}$, and $\frac{1}{r}=\frac{1-\theta}{r_1}+\frac{\theta}{r_2}$.
\end{proposition}

We omit the proof of the proposition which is analogous to \text{Proposition 3.3} in \cite{FMV1}.

\section[Technical lemmas]{Technical lemmas}
To obtain Strichartz estimates, it involves estimating the decay in time on the Schr\"{o}dinger operator. By Fourier inversion formula \eqref{Fourier-inversion}, we may write such operators explicitly into a sum of a list of oscillatory integrals. In order to estimate the oscillatory integrals, we shall recall the stationary phase lemma that will be the central argument.
\begin{lemma} (see \cite{S}, Corollary, p.334)\label{phase} Let $g \in C^\infty([a,b])$ be real-valued such that
\begin{equation*}
|g''(x)|\geq \delta,
\end{equation*}
for any $x\in[a,b]$ with $\delta >0$. Then for any function $h \in C^\infty([a,b])$, there exists a constant $C$ which does not depend on $\delta, a, b, g$ or $h$, such that
\begin{equation*}
\left|\int_a^b e^{ig(x)}h(x)\,dx\right|\leq C\delta^{-1/2}\left[||h||_\infty+\int_a^b|h'(x)|\,dx\right].
\end{equation*}
\end{lemma}
In order to prove the sharpness of the time decay in Theorem \ref{dispersive inequ}, we describe the asymptotic expansion of oscillating integrals.
\begin{lemma}(see \cite{S}, Proposition 6, p.344)\label{asymptotic} Suppose $\phi$ is a smooth function on $\mathbb{R}^p$ and has a nondegenerate critical point at $x_0$. If $\psi$ is supported in a sufficiently small neighborhood of $x_0$, then
\begin{equation*}
\left|\int_{\mathbb{R}^p}e^{it\phi(x)}\psi(x)\,dx\right| \sim |t|^{-p/2}, \text{ as t}\rightarrow \infty.\\
\end{equation*}
\end{lemma}

Next, we will need some estimates of the Laguerre functions:
\begin{lemma}(see \cite{DH}, Lemma 3.2)\label{Laguerre}\indent
\begin{equation}
\left|(\tau \frac{d}{d\tau})^k \mathfrak{L}_m^{(d-1)}(\tau)\right|\leq C_{k,d}(2m+d)^{d-1/4},
\end{equation}
for all $0\leq k\leq d$.
\end{lemma}

\begin{remark} In fact, for $0\leq k\leq d-1$, we have a better estimate
\begin{equation*}
\left|(\tau \frac{d}{d\tau})^k \mathfrak{L}_m^{(d-1)}(\tau)\right|\leq C_{k,d}(2m+d)^{d-1}.\\
\end{equation*}
\end{remark}

Furthermore, we will exploit the following estimates, which can be easily proved by comparing the sums with the corresponding integrals:
\begin{lemma}\label{sum}
Fix $\beta \in \mathbb{R}$. There exists $C_\beta>0$ such that for $A>0$ and $d\in\mathbb{Z}_+$, we have
\begin{align}
\underset{2m+d\geq A}{\underset{m\in\mathbb{N}}{\sum}}(2m+d)^\beta &\leq C_\beta A^{\beta+1}, \quad \beta<-1;\label{sum1}\\
\underset{2m+d\leq A}{\underset{m\in\mathbb{N}}{\sum}}(2m+d)^\beta &\leq C_\beta A^{\beta+1}, \quad \beta>-1.\label{sum2}
\end{align}
\end{lemma}

Finally, we introduce the following properties of the Fourier transform of surface-carried measures.
\begin{theorem}(see \cite{So}, Theorem 1.2.1) \label{measure} Let $S$ be a smooth hypersurface in $\mathbb{R}^p$ with non-vanishing Gaussian curvature and $d\mu$ a $C_0^\infty$ measure on $S$. Suppose that $\Gamma\subset\mathbb{R}^p\setminus\{0\}$ is the cone consisting of all $\xi$ which are normal of some point $x\in S$ belonging to a fixed relatively compact neighborhood $\mathcal{N}$ of \text{supp} $d\mu$. Then,
\begin{align*}
\left(\frac{\partial}{\partial\xi}\right)^\alpha\widehat{d\mu}(\xi)&=O\left((1+|\xi|)^{-M}\right), \quad\forall M\in\mathbb{N},\text{ if }\xi\not\in\Gamma,\\
\widehat{d\mu}(\xi)&=\sum e^{-i(x_j,\xi)}a_j(\xi),\text{ if }\xi\in\Gamma,
\end{align*}
where the (finite) sum is taken over all points $x\in\mathcal{N}$ having $\xi$ as a normal and
\begin{equation*}
\left|\left(\frac{\partial}{\partial\xi}\right)^\alpha a_j(\xi)\right|\leq C_\alpha(1+|\xi|)^{-(p-1)/2-|\alpha|}.
\end{equation*}
\end{theorem}

We shall need the following properties of the Fourier transform of the measure $d\sigma$ on the unit sphere $\mathbb{S}^{p-1}$. Obviously, $\widehat{d\sigma}$ is radial. By Theorem \ref{measure}, we have the radical decay properties of the Fourier transform of the spherical measure.
\begin{lemma}\label{sphere} \indent For any $\xi\in\mathbb{R}^p$, the estimate holds
\begin{equation*}
\widehat{d\sigma}(\xi)=e^{i|\xi|}\phi_+(|\xi|)+e^{-i|\xi|}\phi_-(|\xi|),
\end{equation*}
where
\begin{equation*}
|\phi^{(k)}_\pm(r)|\leq c_k(1+r)^{-(p-1)/2-k}, \text{ for all }r>0,\,k\in\mathbb{N}.
\end{equation*}
\end{lemma}

\section[Dispersive estimates]{Dispersive estimates}
First we prove the weak time dispersion which concerns the Littlewood-Paley functions $\psi_j$ defined in Section 3.
\begin{proposition}\label{Weak-time} There exists a constant $C>0$, which depends only on $d$ and $p$, such that for any $\rho \in [n-1,N-2]$, $j\in\mathbb{Z}$ and $t \in \mathbb{R}^*=\mathbb{R}\backslash\{0\}$ we have:
\begin{equation*}
||e^{it\mathcal{L}}\psi_j||_{L^\infty(G)} \leq C|t|^{-1/2}2^{j\rho}.
\end{equation*}
\end{proposition}

{\bf Proof.} \indent Fixing $t \in \mathbb{R}^*$, $j\in \mathbb{Z}$ and $(z,s)\in G$, by the Fourier inversion formula, we have
\begin{equation*}
\begin{split}
e^{it\mathcal{L}}\psi_j(z,s)&=(\frac{1}{2\pi})^{d+p}\underset{m\in\mathbb{N}}{\sum}\int_{\mathbb{R}^p} e^{-i\lambda s} e^{it\left((2m+d)|\lambda|+|\lambda|^2\right)}\\
&\quad \quad \quad \quad \quad \times R\left(2^{-2j}\left((2m+d)|\lambda|+|\lambda|^2\right)\right)\mathfrak{L}_m^{(d-1)}(\frac{|\lambda|}{2}|z|^2)|\lambda|^d \,d\lambda\\
&=(\frac{1}{2\pi})^{d+p}\underset{m\in\mathbb{N}}{\sum}I_m,
\end{split}
\end{equation*}
where
\begin{equation}
I_m=\int_{\mathbb{R}^p} e^{-i\lambda s} e^{it\left((2m+d)|\lambda|+|\lambda|^2\right)}R\left(2^{-2j}\left((2m+d)|\lambda|+|\lambda|^2\right)\right) \mathfrak{L}_m^{(d-1)}(\frac{|\lambda|}{2}|z|^2)|\lambda|^d \,d\lambda,\label{main-I-m}
\end{equation}
and our assertion simply reads
\begin{equation*}
\underset{m\in \mathbb{N}}{\sum}|I_m|\lesssim \bigg\{ \begin{array}{ll}
|t|^{-1/2}2^{j(2d+p-1)}, &j>0\\
|t|^{-1/2}2^{j(2d+2p-2)},&j\leq0.
\end{array}
\end{equation*}
Putting $\delta=\frac{s}{t}$ and $M=2m+d$, we first integrate \eqref{main-I-m} on $\mathbb{R}^+$, then
\begin{equation*}
\begin{split}
I_m&=\int_{\mathbb{R}^p} e^{it(M|\lambda|+|\lambda|^2-\lambda\cdot\delta)}R\left(2^{-2j}(M|\lambda|+|\lambda|^2)\right) \mathfrak{L}_m^{(d-1)}(\frac{|\lambda|}{2}|z|^2)|\lambda|^d \,d\lambda\\
   &=\int_{S^{p-1}}I_{\epsilon,m}\,d\sigma(\epsilon),
\end{split}
\end{equation*}
where
\begin{equation*}
I_{\epsilon,m}=\int_0^{+\infty} e^{it(M\lambda+\lambda^2-\lambda\epsilon\cdot\delta)}R\left(2^{-2j}(M\lambda+\lambda^2)\right) \mathfrak{L}_m^{(d-1)}(\frac{\lambda}{2}|z|^2)\lambda^{d+p-1}\,d\lambda.
\end{equation*}
Performing the change of variable $x=2^{-2j}M\lambda$, we obtain
\begin{equation*}
I_{\epsilon,m}=2^{j(2d+2p)}K_{\epsilon,m},
\end{equation*}
where
\begin{equation} \label{epsilonK}
K_{\epsilon,m}=\int_0^{+\infty} e^{it2^{2j}G_{j,\sigma,\epsilon,m}(x)} h_{j,z,m}(x)\,dx.
\end{equation}
Here,
\begin{align*}
G_{j,\delta,\epsilon,m}(x)&=x+\frac{2^{2j}}{M^2}x^2-\frac{\epsilon\cdot\delta}{M}x,\\
h_{j,z,m}(x)&=R(x+\frac{2^{2j}}{M^2}x^2)\mathfrak{L}^{(d-1)}_m(\frac{2^{2j-1}x|z|^2}{M})\frac{x^{d+p-1}}{M^{d+p}}.
\end{align*}
So
\begin{equation*}
\text{supp }h_{j,z,m}\subseteq \{x \in \mathbb{R}^+:\frac{1}{2}\leq x+\frac{2^{2j}}{M^2}x^2\leq 4\}=[a_{j,m},b_{j,m}],
\end{equation*}
where
\begin{equation*}
a_{j,m}=\frac{1}{1+\sqrt{1+2^{2j+1}M^{-2}}},\quad b_{j,m}=\frac{8}{1+\sqrt{1+2^{2j+4}M^{-2}}}.
\end{equation*}
Note that
\begin{equation}
a_{j,m}, b_{j,m}\sim \min(1,2^{-j}M) \label{ab}.
\end{equation}
For $x\in [a_{j,m}, b_{j,m}]$, we have
\begin{equation}
G''_{j,\delta,\epsilon,m}(x)=\frac{2^{2j+1}}{M^2}.
\end{equation}
Moreover, by \text{Lemma} \ref{Laguerre} and \eqref{ab}, one can easily verify that
\begin{equation*}
||h_{j,z,m}||_{L^\infty[a_{j,m},b_{j,m}]}+||h'_{j,z,m}||_{L^1[a_{j,m},b_{j,m}]}\lesssim \bigg\{ \begin{array}{ll}
M^{-(p+1)},& M\geq 2^j, \\
2^{-j(d+p-1)}M^{d-2},& M<2^j.
\end{array}
\end{equation*}
Applying Lemma \ref{phase}, we obtain a consistent estimate
\begin{equation*}
 \left|K_{\epsilon,m}\right|\lesssim \bigg\{ \begin{array}{ll}
|t|^{-1/2}2^{-2j}M^{-p},& M\geq 2^j, \\
|t|^{-1/2}2^{-j(d+p+1)}M^{d-1},& M<2^j.
\end{array}
\end{equation*}
Hence, we have
\begin{equation}\label{I_m}
|I_m|\lesssim \bigg\{ \begin{array}{ll}
|t|^{-1/2}2^{j(2d+2p-2)}M^{-p},& M\geq 2^j, \\
|t|^{-1/2}2^{j(d+p+1)}M^{d-1},& M<2^j.
\end{array}
\end{equation}
Noting that $p>1$, for $j\leq 0$, $\underset{m\in \mathbb{N}}{\sum}|I_m|\lesssim |t|^{-1/2}2^{j(2d+2p-2)}$. For $j>0$, $\underset{m\in \mathbb{N}}{\sum}|I_m|\lesssim |t|^{-1/2}2^{j(2d+p-1)}$ follows from \eqref{I_m} by applying \text{Lemma} \ref{sum} separately to the sums $\sum\nolimits_{M\geq 2^j}|I_m|$ and $\sum\nolimits_{M<2^j}|I_m|$.\\\\

\begin{remark}
If we integrate \eqref{main-I-m} first over $S^{p-1}$, we have
\begin{equation*}
I_m=\int_0^{+\infty}\widehat{d\sigma}(\lambda|s|) e^{it(M\lambda+\lambda^2)}R(2^{-2j}(M\lambda+\lambda^2)) \mathfrak{L}_m^{(d-1)}(\frac{\lambda}{2}|z|^2)\lambda^{d+p-1}\,d\lambda.
\end{equation*}
Performing the change of variable $x=2^{-2j}M\lambda$, we obtain
\begin{equation}\label{im}
I_m=2^{j(2d+2p)}\int_0^{+\infty}e^{it2^{2j}(x+\frac{2^{2j}}{M^2}x^2)}\widehat{d\sigma}(\frac{2^{2j}|s|x}{M})h_{j,z,m}(x)\,dx.\\
\end{equation}
It follows from Lemma \ref{sphere} that
\begin{equation}\label{Isum}
I_m=(2\pi)^{p/2}\underset{\pm}{\sum}I_m^\pm,
\end{equation}
where
\begin{equation*}
I_m^\pm=2^{j(2d+2p)}\int_0^{+\infty}e^{it2^{2j}\left(x+\frac{2^{2j}}{M^2}x^2\pm\frac{|s|}{Mt}x\right)}\phi_\pm\left(\frac{2^{2j}|s|}{M}x\right)h_{j,z,m}(x)\,dx.
\end{equation*}
From \text{Lemma} \ref{Laguerre}, Lemma \ref{sphere} and \eqref{ab} one can verify that
\begin{align*}
&\left|\left|\phi_\pm\Big(\frac{2^{2j}|s|}{M}\cdot\Big) h_{j,z,m}\right|\right|_{L^\infty[a_{j,m},b_{j,m}]}+\left|\left|\frac{\partial}{\partial x}(\phi_\pm\Big(\frac{2^{2j}|s|}{M}\cdot) h_{j,z,m}\Big)\right|\right|_{L^1[a_{j,m},b_{j,m}]}\\
&\lesssim \bigg\{ \begin{array}{ll}
|s|^{-(p-1)/2}2^{-j(p-1)}M^{-(p+3)/2},& M\geq 2^j, \\
|s|^{-(p-1)/2}2^{-j(d+3(p-1)/2)}M^{d-1},& M<2^j.
\end{array}
\end{align*}
Exploiting Lemma \ref{phase}, we obtain
\begin{equation} \label{pm}
|I^\pm_m|\lesssim \bigg\{ \begin{array}{ll}
|t|^{-1/2}|s|^{-(p-1)/2}2^{j(2d+p-1)}M^{-(p+1)/2},& M\geq 2^j, \\
|t|^{-1/2}|s|^{-(p-1)/2}2^{j(d+(p-1)/2)}M^{d-1},& M<2^j.
\end{array}
\end{equation}
\end{remark}

To improve the time decay, we will try to apply $p$ times a non-critical phase estimate. We will exploit the following estimates for the derivatives of $h_{j,z,m}$.
\begin{lemma}\label{h} For any $x\in[a_{j,m}, b_{j,m}]$, $0\leq l\leq d$, we have
\begin{equation*}
|h^{(l)}_{j,z,m}(x)|\lesssim \bigg\{ \begin{array}{ll}
M^{-(p+\theta_l)},& M\geq 2^j\\
2^{-j(d+p-l-1)}M^{d-l-\theta_l-1},& M<2^j,\\
\end{array}
\end{equation*}
where $\theta_l=\bigg\{ \begin{array}{ll}
1,& 0\leq l\leq d-1 \\
1/4,& l=d\\
\end{array}$.
\end{lemma}

{\bf Proof.} Recall that
\begin{equation*}
h_{j,z,m}(x)=R(x+\frac{2^{2j}}{M^2}x^2)\mathfrak{L}^{(d-1)}_m(\frac{2^{2j-1}x|z|^2}{M^2})\frac{x^{d+p-1}}{M^{d+p}}.
\end{equation*}
By an induction we get
\begin{align*}
h^{(l)}_{j,z,m}(x)&=\underset{\alpha\in\mathcal{F}}{\sum}A(l,\alpha)R^{(\alpha_1)}(x+\frac{2^{2j}}{M^2}x^2)(1+\frac{2^{2j+1}}{M^2}x)^{\alpha_2}(\frac{2^{2j+1}}{M^2})^{\alpha_3}\\
&\quad \quad \quad \times\left[\left(x\frac{d}{dx}\right)^{\alpha_4}\mathfrak{L}^{(d-1)}_m\right](\frac{2^{2j-1}x|z|^2}{M^2})\frac{x^{d+p-\alpha_5-1}}{M^{d+p}}.
\end{align*}
where $\mathcal{F}=\{\alpha=(\alpha_1, \cdots, \alpha_5)\in \mathbb{N}^5: \alpha_1=\alpha_2+\alpha_3, \alpha_1+\alpha_3+\alpha_5=l,\alpha_4\leq\alpha_5\}$.\\
Applying \text{Lemma} \ref{Laguerre} and \eqref{ab}, \text{Lemma} \ref{h} comes out easily.\\

Now we can prove the sharp time dispersion, which is a big improvement on the time decay in Proposition \ref{Weak-time}.
\begin{proposition}\label{Sharp time dispersion} There exists a constant $C>0$, which depends only on $d$ and $p$, such that for any $j\in\mathbb{Z}$ and $t \in \mathbb{R}^*$ we have:
\begin{equation*}
||e^{it\mathcal{L}}\psi_j||_{L^\infty(G)}\leq C|t|^{-p/2}2^{j(n-1)}.
\end{equation*}
\end{proposition}

{\bf Proof.} \indent From Proposition \ref{Weak-time}, it suffices to prove the case $|t|>1$. Without loss of generality, we can assume that $t>1$.\\
For $j>0$, recall from \eqref{epsilonK} that
\begin{equation*}
K_{\epsilon,m}=\int_0^{+\infty} e^{it2^{2j}G_{j,\delta,\epsilon,m}(x)} h_{j,z,m}(x)\,dx,
\end{equation*}
where
\begin{equation}
G'_{j,\delta,\epsilon,m}(x)=1+\frac{2^{2j+1}}{M^2}x-\frac{\epsilon\cdot\delta}{M}\label{Gfirst},
\end{equation}
and
\begin{equation}
G''_{j,\delta,\epsilon,m}(x)=\frac{2^{2j+1}}{M^2}. \label{Gsecond}
\end{equation}
We divide $\mathbb{N}$ into three (possible empty) disjoint subsets:
\begin{align*}
A_1&=\{m\in\mathbb{N}: M\geq2^j, |\delta|\lesssim M\},\\
A_2&=\{m\in\mathbb{N}: M\geq2^j, |\delta|\gtrsim M\},\\
A_3&=\{m\in\mathbb{N}: M<2^j\}.
\end{align*}
Then our assertion reads:
\begin{equation*}
\underset{m\in A_l}{\sum}|I_m|\lesssim t^{-p/2}2^{j(2d+p-1)},\quad l=1,2,3.
\end{equation*}
For $l=1$, by \eqref{Gfirst}, we obtain
\begin{equation}\label{Dfirst}
|G'_{j,\delta,\epsilon,m}(x)|\gtrsim 1, \text{ for any x} \in [a_{j,m}, b_{j,m}].
\end{equation}
The phase function $G'_{j,\delta,\epsilon,m}(x)$ for $K_{\epsilon,m}$ has no critical points on $[a_{j,m}, b_{j,m}]$. By $Q-$fold integration by parts, we get
\begin{equation}\label{integration}
K_{\epsilon,m}=(it2^{2j})^{-Q}\int_0^{+\infty}e^{it2^{2j}G_{j,\delta,\epsilon,m}(x)}D^Qh_{j,z,m}(x)\,dx,
\end{equation}
where the differential operator $D$ is defined by
\begin{equation*}
Dh_{j,z,m}(x)=\frac{d}{dx}\left(\frac{h_{j,z,m}(x)}{G'_{j,\delta,\epsilon,m}(x)}\right).
\end{equation*}
By a direct induction, we have
\begin{equation*}
D^Qh_{j,z,m}=\underset{k=Q}{\overset{2Q}{\sum}}\underset{\alpha_1+2\alpha_2=k}{\sum}C(\alpha, k, Q)\frac{h^{(\alpha_1)}_{j,z,m}(G''_{j,\delta,\epsilon,m})^{\alpha_2}}{(G'_{j,\delta,\epsilon,m})^k},
\end{equation*}
with $\alpha=(\alpha_1,\alpha_2)\in\{0, 1 , \cdots, Q\}\times \mathbb{N}$.\\
The estimates \eqref{Gsecond} and \eqref{Dfirst} yield
\begin{equation*}
||D^Qh_{j,z,m}||_\infty\lesssim \underset{0\leq \alpha_1\leq Q}{\text{sup}}||h^{(\alpha_1)}_{j,z,m}||_\infty.
\end{equation*}
Applying \text{Lemma} \ref{h}, for all $0\leq Q\leq d$, we obtain
\begin{equation}\label{D1}
||D^Qh_{j,z,m}||_\infty\lesssim \underset{0\leq \alpha_1\leq Q}{\text{sup}}||h^{(\alpha_1)}_{j,z,m}||_\infty \lesssim M^{-(p+1/4)},
\end{equation}
which yields the trivial estimate
\begin{equation}\label{trivial1}
|K_{\epsilon,m}|\lesssim M^{-(p+1/4)}.
\end{equation}
Moreover, it follows from \eqref{ab}, \eqref{integration} and \eqref{D1} that, for any $0\leq Q\leq d$, we get a uniform estimate (with respect to $\epsilon\in S^{p-1}$)
\begin{equation}\label{K1}
|K_{\epsilon,m}|\lesssim {(t2^{2j})}^{-Q}M^{-(p+1/4)}.
\end{equation}
Interpolating \eqref{trivial1} and \eqref{K1}, we get that for all $0\leq\theta\leq d$
\begin{equation*}
|K_{\epsilon,m}|\lesssim {(t2^{2j})}^{-\theta}M^{-(p+1/4)}.
\end{equation*}
Since $p\leq 2d-1$, we have $(p+1)/2\leq d$. Hence, let $\theta=(p+1)/2$,
\begin{equation*}
|K_{\epsilon,m}|\lesssim t^{-p/2}2^{-j(p+1)}M^{-(p+1/4)}.
\end{equation*}
Finally, the desired estimate holds
\begin{equation*}
\underset{m\in A_1}{\sum}|I_m|\lesssim t^{-p/2}2^{j(2d+p-1)}\underset{m\in\mathbb{N}}{\sum}M^{-(p+1/4)}\lesssim t^{-p/2}2^{j(2d+p-1)}.
\end{equation*}
For $l=2$, the estimate \eqref{pm} yields
\begin{equation*}
|I^\pm_m|\lesssim t^{-p/2}2^{j(2d+p-1)}M^{-p}.
\end{equation*}
Then it follows from \eqref{Isum} that
\begin{equation*}
\underset{m\in A_2}{\sum}|I_m|\lesssim t^{-p/2}2^{j(2d+p-1)}\underset{m\in\mathbb{N}}{\sum}M^{-p}\lesssim t^{-p/2}2^{j(2d+p-1)}.
\end{equation*}
For $l=3$, when $|\delta|\gtrsim2^j$, the estimate \eqref{pm} yields
\begin{equation*}
|I_m^{\pm}|\lesssim t^{-p/2}2^{jd}M^{d-1}.
\end{equation*}
Thanks to \eqref{sum2}, we have
\begin{equation*}
\underset{m\in A_3}{\sum}|I_m|\lesssim t^{-p/2}2^{jd}\underset{M<2^j}{\sum}M^{d-1}\lesssim t^{-p/2}2^{2jd}\lesssim t^{-p/2}2^{j(2d+p-1)}.
\end{equation*}
When $|\delta|\lesssim2^j$, the estimate
\begin{equation}\label{Gfirst2}
|G'_{j,\delta,\epsilon,m}(x)|\gtrsim \frac{2^j}{M},
\end{equation}
holds for any $x \in [a_{j,m}, b_{j,m}]$. \\
According to \text{Lemma} \ref{h}, for any $0\leq\alpha_1\leq d$
\begin{equation}\label{hinfty2}
||h^{\alpha_1}_{j,z,m}||_\infty\lesssim2^{-j(d+p-1)}\left(\frac{2^j}{M}\right)^{\alpha_1}M^{d-5/4}.
\end{equation}
It implies the following trivial estimate
\begin{equation}\label{trivial2}
|K_{\epsilon,m}|\lesssim2^{-j(d+p-1)}M^{d-5/4}.
\end{equation}
Furthermore, analogous to the case $r=1$, \eqref{Gsecond}, \eqref{Gfirst2} and \eqref{hinfty2} yield, for any $0\leq Q\leq d$
\begin{align*}
|D^Qh_{j,z,m}|&\lesssim\underset{k=Q}{\overset{2Q}{\sum}}\underset{\alpha_1+2\alpha_2=k}{\sum}\frac{|h^{(\alpha_1)}_{j,z,m}||G''_{j,\delta,\epsilon,m}|^{\alpha_2}}{|G'_{j,\delta,\epsilon,m}|^k}\\
              &\lesssim2^{-j(d+p-1)}M^{d-5/4}.
\end{align*}
We can also obtain a uniform estimate
\begin{equation}\label{K2}
|K_{\epsilon,m}|\lesssim (t2^{2j})^{-Q}2^{-j(d+p-1)}M^{d-5/4}.
\end{equation}
Interpolating \eqref{trivial2} and \eqref{K2}, for any $0\leq\theta\leq d$
\begin{equation*}
|K_{\epsilon,m}|\lesssim {(t2^{2j})}^{-\theta}2^{-j(d+p-1)}M^{d-5/4}.
\end{equation*}
Let $\theta=p/2$,
\begin{equation*}
|K_{\epsilon,m}|\lesssim t^{-p/2}2^{-j(d+2p-1)}M^{d-5/4}.
\end{equation*}
Finally, because of \eqref{sum2} and $p\geq2$, the desired estimate holds
\begin{equation*}
\underset{m\in A_3}{\sum}|I_m|\lesssim t^{-p/2}2^{j(d+1)}\underset{M<2^j}{\sum}M^{d-5/4}\lesssim t^{-p/2}2^{j(2d+3/4)}\lesssim t^{-p/2}2^{j(2d+p-1)}.
\end{equation*}

For $j\leq0$, recall from \eqref{im} that
\begin{equation*}
I_m=2^{j(2d+2p)}\int_0^{+\infty}e^{it2^{2j}\varphi_{j,m}(x)}H_{j,z,s,m}(x)\,dx,
\end{equation*}
where
\begin{align*}
\varphi_{j,m}(x)&=x+\frac{2^{2j}}{M^2}x^2,\\
H_{j,z,s,m}(x)&=\widehat{d\sigma}(\frac{2^{2j}|s|x}{M})h_{j,z,m}(x).
\end{align*}
First we obtain a trivial estimate
\begin{equation}\label{Trivial1}
|I_m|\lesssim2^{j(2d+2p)}M^{-(p+1)}.
\end{equation}
We will discuss it in the following cases.\\

\textsf{Case 1.} $\frac{2^{2j}|s|}{M}\leq1$. In this case, we will exploit the vanishing property of the Fourier transform of the spherical measure at the origin. One can easily get from Lemma \ref{sphere} that for any $k\in\mathbb{N}$,
\begin{equation}
\left|\frac{\partial^k}{\partial x^k}\left(\widehat{d\sigma}(\frac{2^{2j}|s|x}{M})\right)\right|\leq c_k.
\end{equation}
Analogous to $r=1$ when dealing with $j>0$, using integration by parts, for any $0\leq Q\leq d$,
\begin{equation*}
I_m=2^{j(2d+2p)}(it2^{2j})^{-Q}\underset{k=Q}{\overset{2Q}{\sum}}\underset{\alpha_1+2\alpha_2=k}{\sum}C(\alpha, k, Q)\int_0^{+\infty}e^{it2^{2j}\varphi_{j,m}(x)}\frac{H^{(\alpha_1)}_{j,z,s,m}(x)(\varphi''_{j,m}(x))^{\alpha_2}}{(\varphi'_{j,m}(x))^k}\,dx,
\end{equation*}
with $\alpha=(\alpha_1,\alpha_2)\in\{0, 1 , \cdots, Q\}\times \mathbb{N}$.\\
It follows from \eqref{ab}, \text{Lemma} \ref{h} and \eqref{sphere} that
\begin{equation}\label{H}
||H^{(\alpha_1)}_{j,z,s,m}||_\infty\lesssim M^{-(p+1/4)}, \text{ for any }0\leq\alpha_1\leq d.
\end{equation}
Note that
\begin{align}
\varphi'_{j,m}(x)&\gtrsim1, \label{varphi1}\\
\varphi''_{j,m}(x)&=\frac{2^{2j}}{M^2}  \label{varphi2}.
\end{align}
So \eqref{ab}, \eqref{H}, \eqref{varphi1} and \eqref{varphi2} imply that, for any $0\leq Q\leq d$,
\begin{equation}\label{I}
|I_m|\lesssim2^{j(2d+2p)}(t2^{2j})^{-Q} M^{-(p+1/4)}.
\end{equation}
Interpolating \eqref{Trivial1} and \eqref{I}, we get that for any $0\leq\theta\leq d$,
\begin{equation*}
|I_m|\lesssim2^{j(2d+2p)}(t2^{2j})^{-\theta} M^{-(p+1/4)}.
\end{equation*}
Let $\theta=p/2$,
\begin{equation*}
|I_m|\lesssim t^{-p/2}2^{j(2d+p)}M^{-(p+1/4)}.
\end{equation*}

\textsf{Case 2.} $\frac{2^{2j}|s|}{M}>1$. We will use the decay property of the Fourier transform of the spherical measure. It follows from \eqref{Isum} that
\begin{align*}
I_m&=c_p2^{j(2d+2p)}\int_0^{+\infty}e^{it2^{2j}\left(x+\frac{2^{2j}}{M^2}x^2\right)}\left(e^{i\frac{2^{2j}|s|}{M}x}\phi_+\Big(\frac{2^{2j}|s|}{M}x\Big)+e^{-i\frac{2^{2j}|s|}{M}x}\phi_-\Big(\frac{2^{2j}|s|}{M}x\Big)\right)h_{j,z,m}(x)\,dx\\
   &=c_p2^{j(2d+2p)}\int_0^{+\infty}e^{it2^{2j}\left(x+\frac{2^{2j}}{M^2}x^2+\frac{|\delta|}{M}x\right)}\phi_+\Big(\frac{2^{2j}|s|}{M}x\Big)h_{j,z,m}(x)\,dx\\
   &\quad\quad\quad\quad+c_p2^{j(2d+2p)}\int_0^{+\infty}e^{i\frac{t2^{2j}}{M}\left((M-|\delta|)x+\frac{2^{2j}}{M}x^2\right)}\phi_-\Big(\frac{2^{2j}|s|}{M}x\Big)h_{j,z,m}(x)\,dx\\
   &=:B_m^++B_m^-.
\end{align*}
Note from Lemma \ref{sphere} that
\begin{equation}\label{Phi}
\left|\frac{\partial^k}{\partial x^k}\Bigg(\phi_\pm\Big(\frac{2^{2j}|s|}{M}x\Big)\Bigg)\right|\leq c_k\Big(\frac{2^{2j}|s|}{M}\Big)^{-\frac{p-1}{2}}\leq c_k, \text{ for any } k\geq0.
\end{equation}
Analogous to \text{Case 1}, we can get
\begin{equation*}
|B^+_m|\lesssim t^{-p/2}2^{j(2d+p)}M^{-(p/2+1/4)}.
\end{equation*}
For $B_m^-$, let $\Phi_{j,m,-}(x)=(M-|\delta|)x+\frac{2^{2j}}{M}x^2$. Note that if $\Big|M-|\delta|\Big|\gtrsim1$, we have $\left|\Phi'_{j,m,-}(x)\right|\gtrsim1$. Repeating what we have done in \text{Case 1}, it follows from \eqref{ab}, \text{Lemma} \ref{h}, and \eqref{Phi} that
\begin{equation*}
|B^-_m|\lesssim t^{-p/2}2^{j(2d+p)}M^{-(p+1/4)}.
\end{equation*}
For $\Big|M-|\delta|\Big|\lesssim1$, we have $M\thicksim \frac{|s|}{t}$. According to \text{Lemma} \ref{Laguerre}, \eqref{ab} and \eqref{Phi}, Lemma \ref{phase} implies
\begin{align*}
|B^-_m|&\lesssim 2^{j(2d+2p)}\left(\frac{t2^{4j}}{M^2}\right)^{-\frac{1}{2}}\Big(\frac{2^{2j}|s|}{M}\Big)^{-\frac{p-1}{2}}M^{-(p+1)}\\
       &\lesssim 2^{j(2d+2p)}\left(\frac{t2^{4j}}{M^2}\right)^{-\frac{1}{2}}(2^{2j}t)^{-\frac{p-1}{2}}M^{-(p+1)}\\
       &\lesssim t^{-p/2}2^{j(2d+p-1)}M^{-p}.
\end{align*}
Combining \text{Case 1} and \text{Case 2}, for $j\leq0$ and noting $p\geq2$, the estimate holds
\begin{equation*}
|I_m|\lesssim t^{-p/2}2^{j(2d+p-1)}M^{-(p/2+1/4)}.
\end{equation*}
By summing over $m\in\mathbb{N}$ the proposition is proved.\\

From Proposition \ref{Sharp time dispersion}, it is easy to obtain our sharp dispersive inequality (see the proof of Corollary 10 in \cite{FV}).
\begin{corollary}\label{dispersion} There exists a $C>0$, which depends only on $d$ and $p$, such that for any $t \in \mathbb{R}^*$
\begin{equation*}
 ||e^{it\mathcal{L}}u_0||_{L^\infty(G)}\leq C|t|^{-p/2}||u_0||_{\dot{B}^{n-1}_{1,1}}.
\end{equation*}
\end{corollary}

{\bf Proof of Theorem \ref{dispersive inequ}: } The dispersive inequality in Theorem \ref{dispersive inequ} is a direct consequence of Corollary \ref{dispersion}. It suffices to show the sharpness of the time decay. Let $Q\in C_0^\infty(D_0)$ with $Q(d)=1$, where $D_0$ is a small neighborhood of $d$ such that $0\not \in D_0$. Then
\begin{equation*}
\hat{u}_0(\lambda, m)=\begin{cases}
Q(|\lambda|),& m=0\\
0,&m\not=0
\end{cases}
\end{equation*}
determines a solution of the Cauchy problem \eqref{schrodinger} with $f=0$,
\begin{equation*}
u((z,s),t)=e^{it\mathcal{L}}u_0=C\int_{\mathbb{R}^p}e^{-i\lambda\cdot s+it(d|\lambda|+|\lambda|^2) -|\lambda||z|^2/4}Q(|\lambda|)|\lambda|^d\,d\lambda.
\end{equation*}
Consider $u((0,ts_0),t)$ for a fixed $s_0$ such that $|s_0|=3d$.
\begin{equation*}
u((0,ts_0),t)=C\int_{\mathbb{R}^p}e^{it(d|\lambda|+|\lambda|^2-\lambda\cdot s_0)}Q(|\lambda|)|\lambda|^d\,d\lambda.
\end{equation*}
This oscillating integral has a phase function $\phi(\lambda):=d|\lambda|+|\lambda|^2-\lambda\cdot s_0$ with a unique critical point $\lambda_0=\frac{s_0}{3}$ which is not degenerate. Indeed, the Hessian is equal to
\begin{equation*}
H(\lambda)=\frac{1}{|\lambda|}\left((2|\lambda|+d)\delta_{k,l}-\frac{\lambda_k\lambda_l}{|\lambda|^2}d\right)_{1\leq k,l\leq p}.
\end{equation*}
Let $s_0=(0,\cdots,0,3d)$, so $\lambda_0=\frac{s_0}{3}=(0, \cdots, 0, d)$. The Hessian at $\lambda_0$ is
\begin{equation*}
H(\lambda_0)=
\left\{
\begin{array}{llll}
3 & & & \\
& \ddots& &\\
& & 3 & \\
& & & 2 \\
\end{array}
\right\}.
\end{equation*}
Applying Lemma \ref{asymptotic}, we get
\begin{equation*}
u((0,ts_0),t)\sim |t|^{-p/2}.
\end{equation*}

\section[Strichartz estimates]{Strichartz estimates}
\quad We are now to prove our Strichartz estimates.

\begin{theorem} \label{Strichartz}\indent For $i=1,2$, let $q_i, r_i\in[2,\infty]$ and $\rho_i\in\mathbb{R}$ such that
\begin{align*}
&a)\,\,\frac{2}{q_i}=p(\frac{1}{2}-\frac{1}{r_i}); \\
&b)\,\,\rho_i=-(n-1)(\frac{1}{2}-\frac{1}{r_i}),
\end{align*}
except for $(q_i, r_i, p)=(2,\infty,2)$. Then the following estimates are satisfied:
\begin{align*}
||e^{it\mathcal{L}}u_0||_{L^{q_1}(\mathbb{R},\dot{B}^{\rho_1}_{r_1,2})}&\leq C||u_0||_{L^2(G)},\\
||\int_0^te^{i(t-\tau)\mathcal{L})}f(\tau)\,d\tau||_{L^{q_1}((0,T),\dot{B}^{\rho_1}_{r_1,2})}&\leq C||f||_{L^{{q_2}'}((0,T),\dot{B}^{-\rho_2}_{{r_2}',2})}
\end{align*}
where the constant $C>0$ does not depend on $u_0$, $f$ or $T$.
\end{theorem}

Once we have obtained the estimate in Proposition \ref{Sharp time dispersion}, the proof is classical and a good reference is \cite{GV} or \cite{KT}.  A detailed presentation in this framework is also given by \cite{FV} in the proof of \text{Theorem 11} .

Going back to the Schr\"{o}dinger equation \eqref{schrodinger}, Theorem \ref{solution} is straightforward from Theorem \ref{Strichartz}.

\begin{remark} \indent We compare the results by G. Furioli and A. Veneruso \cite{FV}. Theorem \ref{dispersive inequ} and Theorem \ref{solution}, which are general results on H-type groups with the center dimension $p\geq2$, are also compatible with those on the Heisenberg group. Hence, the results in our paper apply to all the H-type groups.
\end{remark}

{\bf Acknowledgements.} The author is supported by National Natural Science Foundation of China under Grant \#11371036, the Specialized Research Fund for the Doctoral Program of Higher Education of China under Grant \#2012000110059 and the Fundamental Research Funds for the Central Universities under Grant \#3102015ZY068. This work was performed while the author studied as a joint Ph.D. student in the Mathematics Department of Christian-Albrechts-Universit\"{a}t zu Kiel. She thanks Professor Detlef M\"{u}ller for the hospitality of his Department.

\begin{flushleft}
\vspace{0.3cm}\textsc{Manli Song\\School of Natural and Applied
Sciences\\Northwestern Polytechnical University\\Xi'an, Shaanxi 710129\\People's
Republic of China\\}
\vspace{0.3cm}\textsc{School of Mathematical
Sciences\\Peking University\\Beijing 100871\\People's
Republic of China\\}
\vspace{0.3cm}
\emph{E-mail address}: \text{mlsong@nwpu.edu.cn}
\end{flushleft}
\end{document}